\documentclass[11pt]{article}
\usepackage{enumerate}
\usepackage{amssymb,a4wide,latexsym,makeidx,epsfig,fleqn}
\usepackage{amsthm}
\usepackage{amsmath}
\usepackage{enumerate}
\usepackage{graphicx}
\newtheorem{theorem}{Theorem}[section]

\newtheorem{example}[theorem]{Example}
\newtheorem{definition}[theorem]{Definition}
\newtheorem{lemma}[theorem]{Lemma}

\begin{document}
\textwidth 150mm \textheight 225mm
\title{The row left rank of a quaternion unit gain graph in terms of maximum degree   \footnote{This work is supported by the National Natural Science Foundations of China (No. 12371348, 12201258), the Postgraduate Research \& Practice Innovation Program of Jiangsu Normal University (No. 2024XKT1702).}}
\author{{ Yong Lu\footnote{Corresponding author.}, Qi Shen}\\
{\small  School of Mathematics and Statistics, Jiangsu Normal University,}\\ {\small  Xuzhou, Jiangsu 221116,
People's Republic
of China.}\\
{\small E-mails: luyong@jsnu.edu.cn, sq\_jsnu@163.com}}

\date{}
\maketitle
\begin{center}
\begin{minipage}{120mm}
\vskip 0.3cm
\begin{center}
{\small {\bf Abstract}}
\end{center}
{\small Let $\Phi=(G,U(\mathbb{Q}),\varphi)$  be a quaternion unit gain graph (or $U(\mathbb{Q})$-gain graph) of order $n$,
$A(\Phi)$ be the adjacency matrix of $\Phi$ and  $r(\Phi)$ be the row left rank of $\Phi$. Let $\Delta$ be the maximum degree of $\Phi$.
In this paper, we prove that $r(\Phi)\geq\frac{n}{\Delta}$. Moreover,  if $\Phi$ is connected, we obtain that $r(\Phi)\geq\frac{n-2}{\Delta-1}$.
All the corresponding extremal graphs are characterized.

 }
\end{minipage}
\end{center}

\section{Introduction }
Let $G=(V(G),E(G))$ be a simple graph, where $V(G)=\{v_{1},v_{2},\ldots,v_{n}\}$ and $E(G)$ are the vertex set and the edge set of $G$, respectively.
The \emph{adjacency matrix} $A(G)$ of $G$ is the symmetric $n\times n$ matrix with entries $a_{ij}=1$ if $v_{i}$ is adjacent to $v_{j}$ and $a_{ij}=0$ otherwise.
The \emph{rank} (resp., \emph{nullity}) of $G$ is the rank (resp., nullity) of $A(G)$, denoted by $r(\Phi)$ (resp., $\eta(G)$).

In $1957$, Collatz et al. \cite{CS} first wanted to characterize all graphs of order $n$ with $\eta(G)>0$.
This question has strong chemical background, because $\eta(G)=0$ is a necessary condition for a so-called conjugated molecule to be chemically stable, where $G$ is the graph representing the carbon-atom skeleton of this molecule.
Up to the present day, the problem has not been solved.
In recent years, researchers have focused on the bounds for the rank (or nullity) of a simple graph in terms of some known graph parameters, such as:
the matching number  (see \cite{FHLL, LGUO, mf, RCZ, SST, WW}); the number of pendant vertices (see \cite{ccz, ctlz, MWTDAM,  wfg});
the maximum degree (see \cite{CB, clt, sl, WGT, WGUO, ZWS}); the girth (see \cite{candl, zwt}) and so on.

Let $e_{ij}$ be the \emph{oriented edge} from $v_{i}$ to $v_{j}$, $\overrightarrow{E}(G)$ be the set of oriented edges and $e_{ij},e_{ji}\in\overrightarrow{E}(G)$.
Even though $e_{ij}$ stands for an edge and an oriented edge simultaneously, it will always be clear in the content.
A \emph{gain graph} is a triple $\Phi=(G,\Omega,\varphi)$ consisting of an \emph{underlying graph} $G=(V,E)$, the \emph{gain group} $\Omega$ and the \emph{gain function} $\varphi:\overrightarrow{E}(G)\rightarrow\Omega$ such that $\varphi(e_{ij})=\varphi(e_{ji})^{-1}=\overline{\varphi(e_{ji})}$.
The \emph{adjacency matrix} $A(\Phi)$ of $\Phi$ is the Hermitian $n\times n$ matrix with entries $h_{ij}=\varphi(e_{ij})=\varphi_{v_{i}v_{j}}$ if $e_{ij}\in\overrightarrow{E}(G)$ and $h_{ij}=0$ otherwise.
The \emph{rank} of $\Phi$ is the rank of $A(\Phi)$, denoted by $r(\Phi)$.

If $\Omega=\{1\}$, then $\Phi$ is a simple graph.
If $\Omega=\{1,-1\}$, then $\Phi$ is a signed graph.
If $\Omega=\{T\in \mathbb{C}: |T|=1\}$, then $\Phi$ is a complex unit gain graph.
Recently, researchers have turned to extend the research of the rank (or nullity) of simple graphs to
signed graphs (see \cite{HHL, lwn, LWZ, WS, wlt}) and complex unit gain graphs (see \cite{hhd, hhy, SK, LUWH,landwu, lwx, LWZ1, landy, REFF, xzw, yqt}).

Let $\mathbb{R}$ and $\mathbb{C}$ be the fields of the real numbers and complex numbers, respectively.
Let $\mathbb{Q}$ be a four-dimensional vector space over $\mathbb{R}$ with an ordered basis, denoted by $1$, $i$, $j$, and $k$.
A \emph{real quaternion}, simply called \emph{quaternion}, is a vector $q=x_{0}+x_{1}i+x_{2}j+x_{3}k\in \mathbb{Q},$
where $x_{0},x_{1},x_{2},x_{3}$ are real numbers and $i,j,k$ satisfy the following conditions:
$$i^{2}=j^{2}=k^{2}=-1;$$
$$ij=-ji=k, jk=-kj=i, ki=-ik=j.$$
From \cite{ZF}, we know that if $x,y$ and $z$ are three different quaternions, then $(xy)^{-1}=y^{-1}x^{-1}$ and $(xy)z=x(yz)$
(Note that $xy\neq yx$, in general).

Let $q=x_{0}+x_{1}i+x_{2}j+x_{3}k\in \mathbb{Q}$.
The \emph{conjugate} $\bar{q}$ (or $q^{\ast}$) of $q$ is $\bar{q}=x_{0}-x_{1}i-x_{2}j-x_{3}k$.
The \emph{modulus} of $q$ is $|q|=\sqrt{q\bar{q}}=\sqrt{x_{0}^{2}+x_{1}^{2}+x_{2}^{2}+x_{3}^{2}}$.
If $q\neq 0$, then the \emph{inverse} of $q$ is $q^{-1}=\frac{\bar{q}}{|q|^{2}}$.
The \emph{real part} of $q$ is $Re(q)=x_{0}$.
The \emph{imaginary part} of $q$ is $Im(q)=x_{1}i+x_{2}j+x_{3}k$.
The \emph{row left (right) rank} of a quaternion matrix $A\in \mathbb{Q}^{m\times n}$ is the maximum number of rows of $A$ that are left (right) linearly independent.
The \emph{column left (right) rank} of a quaternion matrix $A\in \mathbb{Q}^{m\times n}$ is the maximum number of columns of $A$ that are left (right) linearly independent.

Note that quaternions do not satisfy the commutative law of multiplication. By the following lemma and example, we know that the row left rank of a quaternion matrix is not necessarily equal to the row right rank (the same as column left rank and column right rank).

\noindent\begin{lemma}\label{le:1.1}\cite{LWL}
The row left rank of a quaternion matrix $A$ equals the column right rank of $A$. The row right rank of a quaternion matrix $A$ equals the column left rank of $A$.
\end{lemma}

When calculating the row left (right) rank of a quaternion matrix by elementary row operations, we only can multiply a non-zero quaternion on the left (right) side of a row and add it to other rows. Similarly, when calculating the column left (right) rank of a quaternion matrix by elementary column operations, we only can multiply a non-zero quaternion on the left (right) side of a column and add it to other columns.

\noindent\begin{example}\label{ex:1.1}\cite{QNZ}
Let
  \begin{align*}
 A^{'} =\left
 (\begin{array}{cccc}
 1 & i\\
 -i & 1\\
\end{array}
 \right).
 \end{align*}
 Then the row left, the row right, the column left and the column right ranks of $A^{'}$ are equal to 1. Consider now the matrix $A$ obtained by multiply $j$ on the left side of 2-th row of $A^{'}$ and add it to 1-th row of $A^{'}$. That is

\begin{align*}
 A =\left
 (\begin{array}{cccc}
 1-ji & j+i\\
 -i & 1\\
\end{array}
 \right).
 \end{align*}

 The row left rank \begin{align*}
 r (A)
 =r \left(
\begin{array}{cccc}
 1-ji+(j+i)i & j+i-(j+i)\\
 -i & 1\\
\end{array}
 \right)
 =r \left(
\begin{array}{cccc}
 0 & 0\\
 -i & 1\\
\end{array}
 \right)=1.
 \end{align*}

 The row right rank
\begin{align*}
 r (A)
 =r \left(
\begin{array}{cccc}
 1-ji+i(j+i) & j+i-(j+i)\\
 -i & 1\\
\end{array}
 \right)
 =r \left(
\begin{array}{cccc}
 2ij & 0\\
 -i & 1\\
\end{array}
 \right)=2.
 \end{align*}

 The column left rank \begin{align*}
 r (A)
 =r \left(
\begin{array}{cccc}
 1-ji+i(j+i) & j+i\\
 -i+i & 1\\
\end{array}
 \right)
 =r \left(
\begin{array}{cccc}
 2ij & j+i\\
 0 & 1\\
\end{array}
 \right)=2.
 \end{align*}

 The column right rank \begin{align*}
 r (A)
 =r \left(
\begin{array}{cccc}
 1-ji+(j+i)i & j+i\\
 -i+i & 1\\
\end{array}
 \right)
 =r \left(
\begin{array}{cccc}
 0 & j+i\\
 0 & 1\\
\end{array}
 \right)=1.
 \end{align*}
\end{example}

From the example above, we can know that when we multiply a non-zero quaternion on left side of a row of a quaternion matrix $A$, the row left rank of $A$ maintains unchanged, however, the row right rank of $A$ may be various.

If the gain group is $U(\mathbb{Q})=\{q\in \mathbb{Q}: |q|=1\}$, then $\Phi=(G,U(\mathbb{Q}),\varphi)$ (or $G^{\varphi}$ for short) is a \emph{quaternion unit gain graph} (or $U(\mathbb{Q})$-gain graph).
In this paper, the \emph{rank} $r(G^{\varphi})$ of a quaternion matrix $A(G^{\varphi})\in \mathbb{Q}^{n\times n}$ is defined to be the row left rank of $A(G^{\varphi})$.
Belardo et al. \cite{BBCR} studied quaternion unit gain graphs and their associated spectral theories.
Zhou and Lu \cite{QNZ} obtained the relationship between the row left rank of a quaternion unit gain graph and the rank of its underlying graph.
Kyrchei et al. \cite{K} provided a combinatorial description of the determinant of the Laplacian matrix of a quaternion unit gain graph.

We use $P^{\varphi}_{n}$ and $C^{\varphi}_{n}$ to denote a $U(\mathbb{Q})$-gain path and a $U(\mathbb{Q})$-gain cycle on $n$ vertices.
Denoted by $K_{|V_{1}|,|V_{2}|}^{\varphi}$ a $U(\mathbb{Q})$-gain complete bipartite graph, where $V_{1},V_{2}\subseteq V(G^{\varphi})$, $V_{1}\cup V_{2}=V(G^{\varphi})$ and $V_{1}\cap V_{2}=\varnothing$.
Let $G^{\varphi}$ be a $U(\mathbb{Q})$-gain graph with vertex set $V(G^{\varphi})$ and $U\subseteq V(G^{\varphi})$, we denote by $G^{\varphi}-U$ the $U(\mathbb{Q})$-gain graph obtained from $G^{\varphi}$ by removing the vertices in $U$ together with all incident edges. When $U=\{x\}$, we write $G^{\varphi}-U$ simply as $G^{\varphi}-x$. Sometimes we use the notation $G^{\varphi}-H^{\varphi}$ instead of $G^{\varphi}-V(H^{\varphi})$ when $H^{\varphi}$ is an induced subgraph of $G^{\varphi}$. If $G^{\varphi}_{1}$ is an induced subgraph of $G^{\varphi}$ and $x$ is a vertex not in $G^{\varphi}_{1}$, we write the subgraph of $G^{\varphi}$ induced by $V(G^{\varphi}_{1})\cup\{x\}$ simply as $G^{\varphi}_{1}+x$.

For a simple graph $G$, if two vertices $x,y\in V(G)$ are adjacent, then we write $x\sim y$.
For a vertex $x\in V(G)$, $N_{G}(x)=\{y\in V(G): y\sim x\}$ is the \emph{neighbor set} of $x$ in $G$ and $d_{G}(x)=|N_{G}(x)|$ is the \emph{degree} of $x$ in $G$.
If $d_{G}(x)=1$, then $x$ is called a \emph{pendant vertex} of $G$.
The \emph{maximum degree} of $G$ is $\Delta=\max\{d_{G}(x):x\in V(G)\}$.
For a $U(\mathbb{Q})$-gain graph $G^{\varphi}$, $N_{G^{\varphi}}(x)=N_{G}(x)$ and $d_{G^{\varphi}}(x)=d_{G}(x)$. The pendant vertex and the maximum degree of $G^{\varphi}$ is defined to be the pendant vertex and the maximum degree of its underlying graph, respectively.

If no two edges of $M\subseteq E(G)$ are adjacent, then $M$ is called a \emph{matching} of $G$.
If a matching $M$ has maximum cardinality among all matching of $G$, then $M$ is a \emph{maximum matching} of $G$.
If every vertex of $G$ is incident with one edge in a matching $M$, then $M$ is a \emph{perfect matching} of $G$.
The \emph{matching number} of $G$, denoted by $m(G)$, is the cardinality of a maximum matching of $G$.
Obviously, a perfect matching is a maximum matching.
For a matching $M$, an \emph{$M$-alternating path} is a path that alternates between edges in $M$ and edges not in $M$.
A tree is a \emph{perfect-matching tree} (or PM-tree for short) if it has a  perfect matching.
For a $U(\mathbb{Q})$-gain graph $G^{\varphi}$, the matching number of $G^{\varphi}$ is defined to be the matching number of its underlying graph.
A $U(\mathbb{Q})$-gain tree is a $U(\mathbb{Q})$-gain PM-tree if its underlying tree is a PM-tree.

Here, we introduce some results about maximum degree $\Delta$. Fiorini et al. \cite{SF} obtained the upper bound of $\eta(T)$ in terms of $\Delta$ of $T$ when $T$ is a tree of order $n$, that is
$$\eta(T)\leq n-2\lceil\frac{n-1}{\Delta}\rceil.$$
Song et al. \cite{YS} obtained the relationship between $\eta(G)$ and $\Delta$ for a reduced bipartite graph $G$ of order $n$, that is
$$\eta(G)\leq n-2-2\ln_{2}\Delta.$$
Zhou et al. \cite{ZWS} obtained the upper bound of $\eta(G)$ in terms of $\Delta$ of $G$,
$$\eta(G)\leq\frac{\Delta-1}{\Delta}n.$$
The graphs which attain the upper bound are also given. If the graph $G$ is assumed to be connected, they expected that
$$\eta(G)\leq\frac{(\Delta-2)n+2}{\Delta-1},$$
and the upper bound is attained if and only if $G$ is a cycle $C_{n}$ with $n$ divisible by $4$ or a complete bipartite graph $K_{\frac{n}{2},\frac{n}{2}}$. In \cite{CB}, \cite{sl} and \cite{WGT}, this conjecture was proved.
Cheng et al. \cite{clt} characterized the connected graph with nullity
$$\eta(G)=\frac{(\Delta-2)n+i}{\Delta-1},~i=0,1,2.$$

For a complex unit gain graph graph, Lu et al. \cite{landwu}  characterized the relationship between rank and maximum degree.
 Samanta \cite{AS} generalized Lu's results.

In Section $2$, we give some lemmas about $U(\mathbb{Q})$-gain graphs.
In Section $3$, we characterize the relationship between $r(G^{\varphi})$ and $\Delta$ for an arbitrarily $U(\mathbb{Q})$-gain graph $G^{\varphi}$ and characterize the corresponding extremal $U(\mathbb{Q})$-gain graphs.
In Section $4$, we get the relationship between $r(G^{\varphi})$ and $\Delta$ for a connected $G^{\varphi}$ and characterize the
corresponding extremal $U(\mathbb{Q})$-gain graphs.

\section{Preliminaries}

In this section, we list some known results. For  $U(\mathbb{Q})$-gain graphs, we have the following lemmas and definition.

\noindent\begin{lemma}\label{le:2.9}\cite{QNZ}
Let $P^{\varphi}_{n}$ be a $U(\mathbb{Q})$-gain path. If $n$ is odd, then $r(P^{\varphi}_{n})=n-1$. If $n$ is even, then $r(P^{\varphi}_{n})=n$.
\end{lemma}

\noindent\begin{definition}\label{de:2.1}\cite{QNZ}
Let $C_{n}^{\varphi} (n\geq3)$ be a $U(\mathbb{Q})$-gain cycle, denote by $$\varphi(C_{n}^{\varphi})=\varphi_{v_{1}v_{2}}\varphi_{v_{2}v_{3}}\cdots \varphi_{v_{n-1}v_{n}}\varphi_{v_{n}v_{1}}.$$ Then $C_{n}^{\varphi}$ is said to be:
\begin{displaymath}
\left\{\
        \begin{array}{ll}
          \rm Type~1,&  \emph{if}~\varphi(C_{n}^{\varphi})=(-1)^{n/2}~\emph{and}~n~\emph{is~even};\\
          \rm Type~2,& \emph{if}~\varphi(C_{n}^{\varphi})\neq(-1)^{n/2}~\emph{and}~n~\emph{is~even};\\
          \rm Type~3,& \emph{if}~Re\left((-1)^{{(n-1)}/{2}}\varphi(C_{n}^{\varphi})\right)\neq 0~\emph{and}~n~\emph{is~odd};\\
          \rm Type~4,& \emph{if}~Re\left((-1)^{{(n-1)}/{2}}\varphi(C_{n}^{\varphi})\right)=0~\emph{and}~n~\emph{is~odd}.
        \end{array}
      \right.
\end{displaymath}
\end{definition}

\noindent\begin{lemma}\label{le:2.2}\cite{QNZ}
Let $C_{n}^{\varphi}$ be a $U(\mathbb{Q})$-gain cycle of order $n$. Then
\begin{align*}
r(C_{n}^{\varphi})=\left\{\begin{array}{ll}
n-2, & if \;C_{n}^{\varphi}\; is \;of\; Type \;1;\\
n,& if \;C_{n}^{\varphi}\; is \;of\; Type \;2 \;or\; 3;\\
n-1,& if \;C_{n}^{\varphi}\; is \;of\; Type \;4.\\
\end{array}\right.
\end{align*}
\end{lemma}

\noindent\begin{lemma}\label{le:2.3}\cite{QNZ}
Let $G^{\varphi}$ be a $U(\mathbb{Q})$-gain graph and $v$ be a vertex of $G^{\varphi}$. Then $r(G^{\varphi})-2\leq r(G^{\varphi}-v)\leq r(G^{\varphi})$. Moreover, if $G^{\varphi}_{1}$ is an induced subgraph of $G^{\varphi}$, then $r(G^{\varphi}_{1})\leq r(G^{\varphi})$.
\end{lemma}

\noindent\begin{lemma}\label{le:2.5}\cite{QNZ}
Let $G^{\varphi}$ be a $U(\mathbb{Q})$-gain graph and $x$ be a pendant vertex of $G^{\varphi}$. If $y$ is adjacent to  $x$ in $G^{\varphi}$, then $r(G^{\varphi})=r(G^{\varphi}-x-y)+2$.
\end{lemma}

\noindent\begin{lemma}\label{le:2.6}\cite{QNZ}
\begin{enumerate}[(a)]
\item Let $G^{\varphi}$ be a $U(\mathbb{Q})$-gain graph and $G^{\varphi}=G_{1}^{\varphi}\cup G_{2}^{\varphi}\cup \cdots \cup G_{t}^{\varphi}$, where $G_{1}^{\varphi}, G_{2}^{\varphi},\ldots, G_{t}^{\varphi}$ are connected components of $G^{\varphi}$. Then $r(G^{\varphi})=\sum^{t}_{i=1}r(G_{i}^{\varphi})$.
\item Let $G^{\varphi}$ be a $U(\mathbb{Q})$-gain graph with $n$ vertices. Then $r(G^{\varphi})=0$ if and only if $G^{\varphi}$ is a $U(\mathbb{Q})$-gain graph without edges.
\end{enumerate}
\end{lemma}

For a simple graph, we list the following two lemmas.

\noindent\begin{lemma}\label{le:2.7}\cite{ZWS}
Let $G$ be a graph of rank $r$, then it has an induced subgraph $H$ and $|V(H)|=r(H)=r$.
\end{lemma}

\noindent\begin{lemma}\label{le:2.10}\cite{CB}
Let $G$ be a connected graph of rank $r$, then it has a connected induced subgraph $H$ and $|V(H)|=r(H)=r$.
\end{lemma}

Using the same methods of Lemmas \ref{le:2.7} and \ref{le:2.10}, we have the following lemma.

\noindent\begin{lemma}\label{le:2.8}
Let $G^{\varphi}$ be a $U(\mathbb{Q})$-gain graph of rank $r$, then it has an induced subgraph $H^{\varphi}$ and $|V(H^{\varphi})|=r(H^{\varphi})=r$. Moreover, if $G^{\varphi}$ is connected, then we can choose $H^{\varphi}$ to be connected.
\end{lemma}

\section{A lower bound on the row left rank of $G^{\varphi}$}
In this section, we will obtain a lower bound of the row left rank of a $U(\mathbb{Q})$-gain graph in terms of the maximum degree.
First, we will prove the following lemma.

\noindent\begin{lemma}\label{le:3.1}
Let $K^{\varphi}_{a,b}~(a,b\geq2)$ be a $U(\mathbb{Q})$-gain complete bipartite graph and $V(K^{\varphi}_{a,b})=V_{1}\cup V_{2}$, $|V_{1}|=a$, $|V_{2}|=b$. Then $r(K^{\varphi}_{a,b})=2$ if and only if all the $C^{\varphi}_{4}$ in $K^{\varphi}_{a,b}$ are of Type 1.
\end{lemma}
\noindent\textbf{Proof.}
\textbf{Necessity:}
Let\begin{align*}
 &A(K^{\varphi}_{a,b})=\left (
 \begin{array}{ccccccc}
 \textbf{0} & A_{1}\\
 A_{1}^{\ast} & \textbf{0}\\
 \end{array}
 \right),\\
  \end{align*}
where $ A_{1}^{\ast}$ is the conjugate transpose of $A_{1}$.
Let $\alpha_{1},\alpha_{2},\ldots,\alpha_{a}$ be the row vectors of $A_{1}$.
Since $r(K^{\varphi}_{a,b})=2$, we have $r(A_{1})=1$.
Without loss of generality, let $\alpha_{i}=k_{i}\alpha_{1}$ and $k_{i}\neq0~(i=2,3,\ldots,a)$.

Let $u_{1},u_{2}\in V_{1}$ and $v_{1},v_{2}\in V_{2}$.
For convenience we assume $\alpha_{1},\alpha_{2}$ be the vectors corresponding to $u_{1},u_{2}$ in $A_{1}$, respectively.
Let $\varphi_{u_{i}v_{j}}$ be the elements in $A_{1}$ corresponding to the edge $u_{i}v_{j}$ $(i,j\in\{1,2\})$.
Since $\alpha_{2}=k_{2}\alpha_{1}$, we have
$$\varphi_{u_{2}v_{1}}=k_{2}\varphi_{u_{1}v_{1}}~and~\varphi_{u_{2}v_{2}}=k_{2}\varphi_{u_{1}v_{2}}.$$
Let $C^{\varphi}_{4}$ be the $U(\mathbb{Q})$-gain cycle with vertices $u_{1},v_{1},u_{2},v_{2}$ in turn.
Then
\begin{align*}
\varphi(C^{\varphi}_{4})&=\varphi_{u_{1}v_{1}}\varphi_{v_{1}u_{2}}\varphi_{u_{2}v_{2}}\varphi_{v_{2}u_{1}}\\
&=\varphi_{u_{1}v_{1}}(k_{2}\varphi_{u_{1}v_{1}})^{-1}(k_{2}\varphi_{u_{1}v_{2}})\varphi_{v_{2}u_{1}}\\
&=\varphi_{u_{1}v_{1}}\varphi_{u_{1}v_{1}}^{-1}k_{2}^{-1}k_{2}\varphi_{u_{1}v_{2}}\varphi_{v_{2}u_{1}}\\
&=(\varphi_{u_{1}v_{1}}\varphi_{u_{1}v_{1}}^{-1})(k_{2}^{-1}k_{2})(\varphi_{u_{1}v_{2}}\varphi_{v_{2}u_{1}})=1.\\
\end{align*}
By Definition \ref{de:2.1}, $C^{\varphi}_{4}$ is of Type $1$.

\textbf{Sufficiency:} Let $A_{1}$, $\alpha_{1},\alpha_{2},\ldots,\alpha_{a}$ be the same described in the proof of ``Necessity",
$V_{1}=\{u_{1},u_{2},\ldots,u_{a}\}$ and $V_{2}=\{v_{1},v_{2},\ldots,v_{b}\}$.
For convenience we assume $\alpha_{i}$ be the vector corresponding to $u_{i}$ in $A_{1}$ $(i=1,2,\ldots,a)$.

The induced subgraph with vertex set $\{u_{1},u_{2},v_{1},v_{2}\}$ is $C^{\varphi}_{4}$.
Since all the $C^{\varphi}_{4}$ in $K^{\varphi}_{a,b}$ are of Type 1, we have $\varphi(C^{\varphi}_{4})=\varphi_{u_{1}v_{1}}\varphi_{v_{1}u_{2}}\varphi_{u_{2}v_{2}}\varphi_{v_{2}u_{1}}=1.$
So
\begin{align*}
\varphi_{v_{1}u_{2}}\varphi_{u_{2}v_{2}}\varphi_{v_{2}u_{1}}&=\varphi_{v_{1}u_{1}}\\
\varphi_{u_{2}v_{2}}\varphi_{v_{2}u_{1}}&=\varphi_{u_{2}v_{1}}\varphi_{v_{1}u_{1}}.
\end{align*}
Let $k_{2}=\varphi_{u_{2}v_{2}}\varphi_{v_{2}u_{1}}=\varphi_{u_{2}v_{1}}\varphi_{v_{1}u_{1}}$, then $\varphi_{u_{2}v_{1}}=k_{2}\varphi_{u_{1}v_{1}}~and~\varphi_{u_{2}v_{2}}=k_{2}\varphi_{u_{1}v_{2}}.$

The induced subgraph with vertex set $\{u_{1},u_{2},v_{1},v_{3}\}$ is $C^{\varphi}_{4}$.
Since all the $C^{\varphi}_{4}$ in $K^{\varphi}_{a,b}$ are of Type 1, we have $\varphi(C^{\varphi}_{4})=\varphi_{u_{1}v_{1}}\varphi_{v_{1}u_{2}}\varphi_{u_{2}v_{3}}\varphi_{v_{3}u_{1}}=1.$
So
\begin{align*}
\varphi_{v_{1}u_{2}}\varphi_{u_{2}v_{3}}\varphi_{v_{3}u_{1}}&=\varphi_{v_{1}u_{1}}\\
\varphi_{u_{2}v_{3}}\varphi_{v_{3}u_{1}}&=\varphi_{u_{2}v_{1}}\varphi_{v_{1}u_{1}}.
\end{align*}
Then $k_{2}=\varphi_{u_{2}v_{3}}\varphi_{v_{3}u_{1}}$, thus $\varphi_{u_{2}v_{3}}=k_{2}\varphi_{u_{1}v_{3}}.$

Since the induced subgraph with vertex set $\{u_{1},u_{2},v_{1},v_{j}\}$ is $C^{\varphi}_{4}$,
we use the same method to get $\varphi_{u_{2}v_{j}}=k_{2}\varphi_{u_{1}v_{j}}~(j=1,2,\ldots,b)$. Thus $\alpha_{2}=k_{2}\alpha_{1}$.

Repeat the above steps, we can get that $\alpha_{i}=k_{i}\alpha_{1}~(i=2,3,\ldots,a)$, i.e., $r(A_{1})=1$.
Thus $r(K^{\varphi}_{a,b})=2$.
\quad $\square$\\

Now, we get a lower bound of the row left rank of a $U(\mathbb{Q})$-gain graph.

\noindent\begin{theorem}\label{th:3.2}
Let $G^{\varphi}$ be a $U(\mathbb{Q})$-gain graph of order $n$. If $G^{\varphi}$ has no isolated vertexs, then
$$r(G^{\varphi})\geq\frac{n}{\Delta}.$$
\end{theorem}
\noindent\textbf{Proof.}
Let $r(G^{\varphi})=r$. Combining $G^{\varphi}$ with no isolated vertexs and Lemma \ref{le:2.6}$(b)$, we get $r>0$.
By Lemma \ref{le:2.8}, $G^{\varphi}$ has an induced subgraph $G^{\varphi}_{1}$ and $r(G^{\varphi}_{1})=|V(G^{\varphi}_{1})|=r$.
Let $G^{\varphi}_{2}=G^{\varphi}-G^{\varphi}_{1}$. Next we will prove that $d_{G^{\varphi}_{1}}(x)\geq1$ for any $x\in V(G^{\varphi}_{2})$.

Suppose to the contrary that there exists a vertex $x\in V(G^{\varphi}_{2})$ with $d_{G^{\varphi}_{1}}(x)=0$. Since $G^{\varphi}$ has no isolated vertexs, there exists a vertex $y\in V(G^{\varphi}_{2})$ and $y$ is adjacent to $x$.
Let $G^{\varphi}_{3}=G^{\varphi}_{1}+x+y$.
By Lemma \ref{le:2.3}, we have $r(G^{\varphi}_{3})\leq r(G^{\varphi}).$
Combining $x$ is a pendant vertex of $G^{\varphi}_{3}$ and Lemma \ref{le:2.5}, we have
$$r(G^{\varphi}_{3})=r(G^{\varphi}_{3}-x-y)+2=r(G^{\varphi}_{1})+2=r+2>r(G^{\varphi}),$$
a contradiction.

Let $E_{1}=\{uv|u\in V(G_{1})~and~v\in V(G_{2})\}$. Then
\begin{equation}
|E_{1}|\geq |V(G^{\varphi}_{2})|=n-r.
\end{equation}
Since $r(G^{\varphi}_{1})=|V(G^{\varphi}_{1})|=r$, we have $d_{G^{\varphi}_{1}}(u)\geq1$ for any $u\in V(G^{\varphi}_{1})$.
Since $d_{G^{\varphi}}(u)\leq\Delta$, we have
\begin{equation}
|E_{1}|=\sum_{u\in V(G^{\varphi}_{1})}(d_{G^{\varphi}}(u)-d_{G^{\varphi}_{1}}(u))\leq r(\Delta-1).
\end{equation}
Combining $(1)$ and $(2)$, we have
$$n-r\leq|E_{1}|\leq r(\Delta-1),$$
so $r(G^{\varphi})=r\geq\frac{n}{\Delta}$.
\quad $\square$\\

In the following, we will characterize the $U(\mathbb{Q})$-gain graph $G^{\varphi}$ with $r(G^{\varphi})=\frac{n}{\Delta}$.

\noindent\begin{theorem}\label{th:3.3}
Let $G^{\varphi}$ be a $U(\mathbb{Q})$-gain graph of order $n$. If $G^{\varphi}$ has no isolated vertexs, then $r(G^{\varphi})=\frac{n}{\Delta}$ if and only if $G^{\varphi}=\frac{n}{2\Delta}K^{\varphi}_{\Delta,\Delta}$, and each $C^{\varphi}_{4}$ (if any) in $K^{\varphi}_{\Delta,\Delta}$ is of Type $1$.
\end{theorem}
\noindent\textbf{Proof.}
\textbf{Sufficiency:} Let $G^{\varphi}=\frac{n}{2\Delta}K^{\varphi}_{\Delta,\Delta}$, and each $C^{\varphi}_{4}$ (if any) in $K^{\varphi}_{\Delta,\Delta}$ be of Type $1$.

If $\Delta=1$, by Lemmas \ref{le:2.6}$(a)$ and \ref{le:2.9}, then $r(G^{\varphi})=\frac{n}{2}r(K^{\varphi}_{1,1})=n=\frac{n}{\Delta}$.

If $\Delta\geq2$, by Lemmas \ref{le:2.6}$(a)$ and \ref{le:3.1}, then $r(G^{\varphi})=\frac{n}{2\Delta}r(K^{\varphi}_{\Delta,\Delta})=\frac{n}{\Delta}$.

\textbf{Necessity:} Let $G^{\varphi}_{1}$, $G^{\varphi}_{2}$ and $E_{1}$ be the same described in Theorem \ref{th:3.2}.
Since $r(G^{\varphi})=\frac{n}{\Delta}$, the inequalities $(1)$ and $(2)$ are became equalities, i.e.,
$|E_{1}|=n-r$, $d_{G^{\varphi}_{1}}(u)=1$ and $d_{G^{\varphi}}(u)=\Delta$ for any $u\in V(G^{\varphi}_{1})$.
So $G^{\varphi}_{1}=\frac{r}{2}K^{\varphi}_{1,1}$.

If $\Delta=1$, then $G^{\varphi}=\frac{n}{2}K^{\varphi}_{1,1}$.

If $\Delta\geq2$, then $d_{G^{\varphi}_{2}}(u)=\Delta-1\geq 1$ for any $u\in V(G^{\varphi}_{1})$.
Let $u_{1}u_{2}\in E(G_{1})$,
$$N_{G^{\varphi}_{2}}(u_{1})=\{v_{1},v_{2},\ldots,v_{\Delta-1}\}~\textrm{and}~N_{G^{\varphi}_{2}}(u_{2})=\{w_{1},w_{2},\ldots,w_{\Delta-1}\}.$$
Since $|E_{1}|=n-r$, we have $v_{i}\neq w_{j}$ for any $1\leq i,j\leq \Delta-1$.
Now we will proved that $v_{i}$ is adjacent to $w_{j}$ for any $1\leq i,j\leq \Delta-1$.

Suppose $v_{i}$ is not adjacent to $w_{j}$. Let $G^{\varphi}_{4}=G^{\varphi}_{1}+v_{i}+w_{j}$.
By Lemma \ref{le:2.3}, we have $r(G^{\varphi}_{4})\leq r(G^{\varphi}).$
By Lemma \ref{le:2.5}, we have
$$r(G^{\varphi}_{4})=r(G^{\varphi}_{4}-u_{1}-u_{2}-v_{i}-w_{j})+4=r(G^{\varphi}_{1}-u_{1}-u_{2})+4=r(G^{\varphi}_{1})+2>r,$$
a contradiction.

Since the maximum degree of $G^{\varphi}$ is $\Delta$, for all $i,j=1,2,\ldots,\Delta-1$, $$N_{G^{\varphi}}(v_{i})=\{u_{1},w_{1},w_{2},\ldots,w_{\Delta-1}\}~\textrm{and}~N_{G^{\varphi}}(w_{j})=\{u_{2},v_{1},v_{2},\ldots,v_{\Delta-1}\}.$$
Then the induced subgraph with vertex set $\{\underbrace{u_{1},w_{1},\ldots,w_{\Delta-1}}_{\Delta~ vertices},\underbrace{u_{2},v_{1},\ldots,v_{\Delta-1}}_{\Delta~vertices}\}$ is $K^{\varphi}_{\Delta,\Delta}$. Hence $G^{\varphi}=\frac{n}{2\Delta}K^{\varphi}_{\Delta,\Delta}$.

Since $r(G^{\varphi})=\frac{n}{\Delta}$, we have $r(K^{\varphi}_{\Delta,\Delta})=2$ for each $K^{\varphi}_{\Delta,\Delta}$. By Lemma \ref{le:3.1}, each $C^{\varphi}_{4}$ in $K^{\varphi}_{\Delta,\Delta}$ is of Type $1$.
\quad $\square$

\section{A lower bound on the row left rank of a connected $G^{\varphi}$}
A natural problem is that: if $G^{\varphi}$ is connected, what is the lower bound of the rank of $G^{\varphi}$?
Firstly, we will introduce some lemmas.

\noindent\begin{lemma}\label{le:4.2}\cite{DM}
Let $T$ be a tree with $n$ vertices, matching number $m$ and rank $r$. Then $r=2m$.
\end{lemma}

\noindent\begin{lemma}\label{le:4.3}\cite{QNZ}
Let $T^{\varphi}$ be a $U(\mathbb{Q})$-gain tree of order $n$. Then $A(T^{\varphi})$ and $A(T)$ have the same rank.
\end{lemma}

By Lemmas \ref{le:4.2} and \ref{le:4.3}, $T^{\varphi}$ is a $U(\mathbb{Q})$-gain PM-tree if and only if $|V(T^{\varphi})|=2m=r(T^{\varphi})$.
Using the same method of Lemma $2.4$ in \cite{CB}, we have the following lemma.

\noindent\begin{lemma}\label{le:4.1}
Let $T^{\varphi}$ be a $U(\mathbb{Q})$-gain PM-tree and $P^{\varphi}=x_{1}x_{2}\ldots x_{k}$ be the longest path of $T^{\varphi}$. If $k\geq3$, then $d_{T^{\varphi}}(x_{1})=d_{T^{\varphi}}(x_{k})=1$ and $d_{T^{\varphi}}(x_{2})=d_{T^{\varphi}}(x_{k-1})=2$.
\end{lemma}

If $G^{\varphi}$ is connected of order $n$ and $\Delta=1$, i.e., $G^{\varphi}=P^{\varphi}_{2}$, then $r(P^{\varphi}_{2})=2$.
Next, we will get a lower bound of the rank of a connected $U(\mathbb{Q})$-gain graph $G^{\varphi}$ with $\Delta\geq2$.

\noindent\begin{theorem}\label{th:4.2}
Let $G^{\varphi}$ be a connected $U(\mathbb{Q})$-gain graph of order $n$ and $\Delta\geq2$. Then
$$r(G^{\varphi})\geq\frac{n-2}{\Delta-1}.$$
\end{theorem}
\noindent\textbf{Proof.}
Let $r(G^{\varphi})=r$. Combining $G^{\varphi}$ is connected and Lemma \ref{le:2.8}, $G^{\varphi}$ has a connected induced subgraph $H^{\varphi}$ and $|V(H^{\varphi})|=r(H^{\varphi})=r$.
Let $K^{\varphi}=G^{\varphi}-H^{\varphi}$. Next we will prove that $d_{H^{\varphi}}(x)\geq1$ for any $x\in V(K^{\varphi})$.

Suppose to the contrary that there exists a vertex $x\in V(K^{\varphi})$ with $d_{H^{\varphi}}(x)=0$. Since $G^{\varphi}$ is connected, there exists a vertex $y\in V(K^{\varphi})$ and $y$ is adjacent to $x$.
Let $H^{\varphi}_{1}=H^{\varphi}+x+y$.
By Lemma \ref{le:2.3}, we have $r(H^{\varphi}_{1})\leq r(G^{\varphi}).$
Combining $x$ is a pendant vertex of $H^{\varphi}_{1}$ and Lemma \ref{le:2.5}, we have
$$r(H^{\varphi}_{1})=r(H^{\varphi}_{1}-x-y)+2=r(H^{\varphi})+2=r+2>r(G^{\varphi}),$$
a contradiction.

Let $E_{1}=\{uv|u\in V(H)~and~v\in V(K)\}$. Then
\begin{equation}
|E_{1}|\geq |V(K^{\varphi})|=n-r.
\end{equation}
Since $d_{G^{\varphi}}(u)\leq\Delta$, we have
\begin{equation}
|E_{1}|=\sum_{u\in V(H^{\varphi})}(d_{G^{\varphi}}(u)-d_{H^{\varphi}}(u))\leq r\Delta-2|E(H)|.
\end{equation}
Since $H^{\varphi}$ is connected, we have $|E(H)|\geq|V(H^{\varphi})|-1$. So
\begin{equation}
|E_{1}|\leq r\Delta-2(|V(H^{\varphi})|-1)=r\Delta-2(r-1)=(\Delta-2)r+2.
\end{equation}
Combining $(3)$ and $(5)$, we have
$$n-r\leq|E_{1}|\leq (\Delta-2)r+2,$$
thus $r(G^{\varphi})=r\geq\frac{n-2}{\Delta-1}$.
\quad $\square$\\

Finally, we will characterize the connected $U(\mathbb{Q})$-gain graph $G^{\varphi}$ with $r(G^{\varphi})=\frac{n-2}{\Delta-1}$.

Let $T^{\varphi}$ be a $U(\mathbb{Q})$-gain PM-tree with perfect matching $M$ and $v$ be a pendant vertex of $T^{\varphi}$.
A pendant vertex $v^{\ast}$ is called a \emph{dual pendant vertex} of $v$ in $T^{\varphi}$ if the path between $v$ and $v^{\ast}$ in $T^{\varphi}$ is an $M$-alternating path.
Since the edge incident to pendant vertex in $T^{\varphi}$ must belong to $M$, every pendant vertex of $T^{\varphi}$ has least one dual pendant vertex.

\noindent\begin{theorem}\label{th:4.3}
Let $G^{\varphi}$ be connected a $U(\mathbb{Q})$-gain graph of order $n$ and $\Delta\geq2$. Then $r(G^{\varphi})=\frac{n-2}{\Delta-1}$ if and only if $G^{\varphi}$ is a $U(\mathbb{Q})$-gain cycle which is of Type $1$, or $G^{\varphi}=K^{\varphi}_{\frac{n}{2},\frac{n}{2}}$, and each $C^{\varphi}_{4}$ in $G^{\varphi}$ is of Type $1$.
\end{theorem}
\noindent\textbf{Proof.}
\textbf{Sufficiency:} If $G^{\varphi}=C^{\varphi}_{n}$ is of Type $1$, by Lemma \ref{le:2.2}, then $r(G^{\varphi})=n-2=\frac{n-2}{\Delta-1}$.

If $G^{\varphi}=K^{\varphi}_{\frac{n}{2},\frac{n}{2}}$, and each $C^{\varphi}_{4}$ in $G^{\varphi}$ is of Type $1$, by Lemma \ref{le:3.1}, then $r(G^{\varphi})=2=\frac{n-2}{\frac{n}{2}-1}=\frac{n-2}{\Delta-1}$.

\textbf{Necessity:} Let $r(G^{\varphi})=r$. Combining $G^{\varphi}$ is connected and Lemma \ref{le:2.6}$(b)$, we have $2\leq r\leq n$.

\textbf{Case 1.} $\Delta=2$.

Since $G^{\varphi}$ is connected, $G^{\varphi}=C^{\varphi}_{n}$ or $G^{\varphi}=P^{\varphi}_{n}$.
By Lemmas \ref{le:2.9} and \ref{le:2.2}, we have $r(G^{\varphi})\geq n-2=\frac{n-2}{\Delta-1}$,
with equality if and only if $G^{\varphi}$ is a $U(\mathbb{Q})$-gain cycle which is of Type $1$.

\textbf{Case 2.} $\Delta\geq3$.

Let $H^{\varphi}$, $K^{\varphi}$ and $E_{1}$ be the same described in Theorem \ref{th:4.2}.
Since $r(G^{\varphi})=\frac{n-2}{\Delta-1}$, the inequalities $(3)$, $(4)$ and $(5)$ are became equalities, i.e.,
\begin{enumerate}[(1)]
\item $|E_{1}|=n-r$;
\item $|E(H)|=|V(H^{\varphi})|-1$, i.e., $H^{\varphi}$ is a $U(\mathbb{Q})$-gain tree;
\item $d_{G^{\varphi}}(x)=\Delta$ for each vertex $x\in V(H^{\varphi})$.
\end{enumerate}

\textbf{Subcase 2.1.} $r=2$. Let $H^{\varphi}=P^{\varphi}_{2}$ and $V(H^{\varphi})=\{u,v\}$.
By $(3)$, $d_{K^{\varphi}}(u)=d_{K^{\varphi}}(v)=\Delta-1$.
Let
$$N_{K^{\varphi}}(u)=\{v_{1},v_{2},\ldots,v_{\Delta-1}\}~\textrm{and}~N_{K^{\varphi}}(v)=\{u_{1},u_{2},\ldots,u_{\Delta-1}\}.$$
Since $(1)$, we have $v_{i}\neq u_{j}$ for any $1\leq i,j\leq \Delta-1$.
Then $v_{i}$ is adjacent to $u_{j}$ for any $1\leq i,j\leq \Delta-1$.
Suppose $v_{i}$ is not adjacent to $u_{j}$. Let $H^{\varphi}_{2}=H^{\varphi}+v_{i}+u_{j}$, then $H^{\varphi}_{2}=P^{\varphi}_{4}$.
By Lemma \ref{le:2.3}, we have $r(H^{\varphi}_{2})\leq r.$
By Lemma \ref{le:2.9}, we have $r(H^{\varphi}_{2})=4>r$, a contradiction.
Since the maximum degree of $G^{\varphi}$ is $\Delta$, for all $i,j=1,2,\ldots,\Delta-1$, $$N_{G^{\varphi}}(v_{i})=\{u,u_{1},u_{2},\ldots,u_{\Delta-1}\}~\textrm{and}~N_{G^{\varphi}}(u_{j})=\{v,v_{1},v_{2},\ldots,v_{\Delta-1}\}.$$
Since $G^{\varphi}$ is connected, there are no vertices in $K^{\varphi}$ except the vertices
$v_{1},\ldots,v_{\Delta-1},u_{1},\ldots,u_{\Delta-1}.$
Then $V(G^{\varphi})=\{\underbrace{u,u_{1},\ldots,u_{\Delta-1}}_{\Delta~vertices},\underbrace{v,v_{1},\ldots,v_{\Delta-1}}_{\Delta~vertices}\}$.
Hence $G^{\varphi}=K^{\varphi}_{\frac{n}{2},\frac{n}{2}}$.
Since $r(G^{\varphi})=2$, so by Lemma \ref{le:3.1}, we have each $C^{\varphi}_{4}$ in $G^{\varphi}$ is of Type $1$.

\textbf{Subcase 2.2.} $r\geq3$. Since $r(H^{\varphi})=|V(H^{\varphi})|=r$, we have $H^{\varphi}$ is a $U(\mathbb{Q})$-gain PM-tree.
Let $P^{\varphi}=x_{1}x_{2}\ldots x_{k}$ be the longest path of $H^{\varphi}$.
Since $|V(H^{\varphi})|=r\geq3$, we have $k\geq3$.
By Lemma \ref{le:4.1}, $d_{H^{\varphi}}(x_{1})=1$ and $d_{H^{\varphi}}(x_{2})=2$.
Let $x^{\ast}$ be a dual pendant vertex of $x_{1}$ in $H^{\varphi}$.
Since $(3)$, we may suppose that
$N_{K^{\varphi}}(x_{1})=\{u_{1},u_{2},\ldots,u_{\Delta-1}\},$
$N_{K^{\varphi}}(x_{2})=\{v_{1},v_{2},\ldots,v_{\Delta-2}\}$
and $N_{K^{\varphi}}(x^{\ast})=\{w_{1},w_{2},\ldots,w_{\Delta-1}\}$.
Since $(1)$, $N_{H^{\varphi}}(u_{i})=\{x_{1}\}$, $N_{H^{\varphi}}(v_{j})=\{x_{2}\}$ and $N_{H^{\varphi}}(w_{i})=\{x^{\ast}\}$
$(i=1,2,\ldots,\Delta-1;~j=1,2,\ldots,\Delta-2)$.

Next we will prove that $u_{1}$ is adjacent to $v_{i}$ for any $1\leq i\leq\Delta-2$.
Suppose $u_{1}$ is not adjacent to $v_{i}$.
Since $H^{\varphi}$ is a $U(\mathbb{Q})$-gain PM-tree, we have $H^{\varphi}-x_{1}+v_{i}$ is a $U(\mathbb{Q})$-gain PM-tree and $r(H^{\varphi}-x_{1}+v_{i})=|V(H^{\varphi}-x_{1}+v_{i})|=r$.
By Lemma \ref{le:2.3}, we have $r(H^{\varphi}+u_{1}+v_{i})\leq r(G^{\varphi}).$
By Lemma \ref{le:2.5}, we have
$$r(H^{\varphi}+u_{1}+v_{i})=r(H^{\varphi}-x_{1}+v_{i})+2=r+2>r,$$
a contradiction.

Finally, we will prove that $u_{1}$ is adjacent to $w_{j}$ for any $1\leq j\leq\Delta-1$.
Suppose $u_{1}$ is not adjacent to $w_{j}$.
Since $H^{\varphi}$ is a $U(\mathbb{Q})$-gain PM-tree, $H^{\varphi}$ has a perfect matching $M$.
Let $H^{\varphi}_{3}=H^{\varphi}-x^{\ast}+u_{1}$. Then $H^{\varphi}_{3}$ is a $U(\mathbb{Q})$-gain tree.
Let $P^{\varphi}_{1}$ be an $M$-alternating path between $x_{1}$ and $x^{\ast}$.
Let $M_{1}=(M\cup E(P^{\varphi}_{1}))\setminus(M\cap E(P^{\varphi}_{1}))$. Then $M_{1}$ is a perfect matching for $H^{\varphi}-x_{1}-x^{\ast}$.
Thus $M_{2}=M_{1}\cup\{x_{1}u_{1}\}$ is a perfect matching for $H^{\varphi}_{3}$.
So $H^{\varphi}_{3}$ is a $U(\mathbb{Q})$-gain PM-tree and $r(H^{\varphi}_{3})=r$.
Let $H^{\varphi}_{4}=H^{\varphi}+u_{1}+w_{j}$.
By Lemma \ref{le:2.3}, we have $r(H^{\varphi}_{4})\leq r(G^{\varphi}).$
By Lemma \ref{le:2.5}, we have
$$r(H^{\varphi}_{4})=r(H^{\varphi}_{4}-x^{\ast}-w_{j})+2=r(H^{\varphi}_{3})+2=r+2>r,$$
a contradiction.

Thus $d_{G^{\varphi}}(u_{1})=1+(\Delta-2)+(\Delta-1)=2\Delta-2>\Delta$ (since $\Delta\geq3$), a contradiction.

So $r=2$, $G^{\varphi}=K^{\varphi}_{\frac{n}{2},\frac{n}{2}}$ and each $C^{\varphi}_{4}$ in $G^{\varphi}$ is of Type $1$.
\quad $\square$

\textbf{Declaration of competing interest}

The authors declare that they have no known competing financial interests or personal relationships that could have appeared to influence the work reported in this paper.

\textbf{Date availability}

No date was used for the research described in the article.


\begin{thebibliography}{99}

\bibitem{BBCR} F. Belardo, M. Brunetti, N.J. Coble, N. Reff, H. Skogman, Spectral of quaternion unit gain graphs, Linear Algebra Appl. 632(2022) 15--49.


\bibitem{ccz} S. Chang, A. Chang, Y.R. Zheng, The leaf-free graphs with nullity $2c(G)-1$, Discrete Appl. Math. 277(2020) 44-54.


\bibitem{DM} D. M. Cvetkovi, I. Gutman, The algebraic multiplicity of the number zero in the spectrum of a bipartite graph, Mat. vesnik. 9(1972) 141-150.

\bibitem{candl} S. Chang, J.X. Li, Graphs $G$ with nullity $n(G)-g(G)-1$, Linear Algebra Appl. 642(2022) 251-263.

\bibitem{CB} B. Cheng, M.H. Liu, B.L. Liu, Proof of a conjecture on the nullity of a connected graph in terms of order and maximum degree, Linear Algebra Appl. 587(2020) 291-301.


\bibitem{clt} B. Cheng, M.H. Liu, B.-S. Tam, On the nullity of a connected graph in terms of order and maximum degree, Linear Algebra Appl. 632(2022) 192-232.

\bibitem{CS} L. Collatz, U. Sinogowitz, Spektren endlicher grafen, Abh. Math. Sem. Univ. Hamburg. 21(1957) 63-77.

\bibitem{ctlz} S. Chang, B.-S. Tam, J.X. Li, Y.R. Zheng, Graphs $G$ with nullity $2c(G)+p(G)-1$,  Discrete Appl. Math. 311(2022) 38-58.






\bibitem{SF} S. Fiorini, I. Gutman, I. Sciriha, Trees with maximum nullity, Linear Algebra Appl. 397(2005) 245-251.

\bibitem{FHLL} Z.M. Feng, J. Huang, S.C. Li, X.B. Luo, Relationship between the rank and the matching number of a graph, Appl. Math. Comput.  354(2019) 411-421.


\bibitem{hhd} S.J. He, R.X. Hao, F.M. Dong, The rank of a complex unit gain graph in terms of the matching number, Linear Algebra Appl. 589(2020) 158-185.

\bibitem{HHL} S.J. He, R.X. Hao, H.-J. Lai, Bounds for the matching number and cyclomatic number of a signed graph in terms of rank, Linear Algebra Appl. 572(2019) 273-291.

\bibitem{hhy} S.J. He, R.X. Hao, A.M. Yu, Bounds for the rank of a complex unit gain graph in terms of the independence number, Linear  Multilinear Algebra 70(2022) 1382-1402.

\bibitem{SK} S. Khan, On connected $\mathbb{T}$-gain graphs with rank equal to girth, Linear Algebra Appl. 688(2024) 232-243.

\bibitem{K} I.I. Kyrchei, E. Treister, V.O. Pelykh, The determinant of the Laplacian matrix of a quaternion unit gain graph, Discrete Math. 347(2024) 113955.

\bibitem{LWL} W.L. Li, Quaternion matrices, National Defense Science and Technology University, 2002.

\bibitem{LGUO} X. Li, J.-M. Guo, No graph with nullity $\eta(G)=|V(G)|-2m(G)+2c(G)-1$, Discrete Appl. Math. 268(2019) 130-136.


\bibitem{LUWH} L. Lu, J.F. Wang, Q.X. Huang, Complex unit gain graphs with exactly one positive eigenvalue,  Linear Algebra Appl. 608(2021)   270-281.



\bibitem{landwu} Y. Lu, J.W. Wu, Bounds for the rank of a complex unit gain graph in terms of its maximum degree, Linear Algebra Appl. 610(2021) 73-85.

\bibitem{lwn} Y. Lu, J.W. Wu, No signed graph with the nullity $\eta(G,\sigma)=|V(G)|-2m(G)+2c(G)-1$, Linear Algebra Appl. 615(2021) 175-193.

\bibitem{lwx} Y. Lu, L.G. Wang, P. Xiao, Complex unit gain bicyclic graphs with rank 2,3 or 4, Linear Algebra Appl. 523(2017) 169-186.

\bibitem{LWZ} Y. Lu, L.G. Wang, Q.N. Zhou, The rank of a signed graph in terms of the rank of its underlying graph, Linear Algebra Appl. 538(2018) 166-186.

\bibitem{LWZ1} Y. Lu, L.G. Wang, Q.N. Zhou, The rank of a complex unit gain graph in terms of the rank of its underlying graph,  J. Comb. Optim. 38(2019) 570-588.






\bibitem{landy} S.C. Li, T. Yang, On the relation between the adjacency rank of a complex unit gain graph and the matching number of its underlying graph, Linear  Multilinear Algebra 70(2022) 1768-1787.

\bibitem{mf} X.B. Ma, X.W. Fang, An improved lower bound for the nullity of a graph in terms of matching number, Linear Multilinear Algebra 68(2020) 1983-1989.

\bibitem{MWTDAM} X.B. Ma,  D. Wong, F.L. Tian, Nullity of a graph in terms of the dimension of cycle space and the number of pendant vertices, Discrete Appl. Math. 215(2016) 171-176.

\bibitem{REFF} N. Reff, Spectral properties of complex unit gain graphs, Linear Algebra Appl. 436(2012) 3165-3176.

\bibitem{RCZ} S. Rula, A. Chang, Y.R. Zheng, The extremal graphs with respect to their nullity, J. Inequal. Appl. 2016(2016) 1-13.

\bibitem{AS} A. Samanta, On bounds of $A_{\alpha}$-eigenvalue multiplicity and the rank of a complex unit gain graph, Discrete Math. 346(2023) 113503.

\bibitem{sl} W.T. Sun, S.C. Li, A short proof of Zhou, Wong and Sun's conjecture, Linear Algebra Appl. 589(2020) 80-84.

\bibitem{SST} Y.Z. Song, X.Q. Song, B.-T. Tam, A Characterization of graphs $G$ with nullity $V(G)-2m(G)+2c(G)$, Linear Algebra Appl. 465(2015) 363-375.

\bibitem{YS} Y. Song, X. Song, C. Zhang, An upper bound for the nullity of a bipartite graph in terms of its maximum degree, Linear Multilinear Algebra 64(2016) 1107-1112.



\bibitem{WS} S.J. Wang, Relation between the rank of a signed graph and the rank of its underlying graph,  Linear Multilinear Algebra 67(2019) 2520-2539.

\bibitem{wfg} L. Wang, X.W. Fang, X.Y. Geng, Graphs with nullity $2c(G)+p(G)-1$, Discrete Math. 345(2022) 112786.

\bibitem{WGT} L. Wang, X. Geng, Proof of a conjecture on the nullity of a graph, J. Graph Theory. 95(2020) 586-593.

\bibitem{WW} L. Wang, D. Wong, Bounds for the matching number, the edge charomatic numbber and the independence number of a graph in terms of rank, Discrete Appl. Math. 166(2014) 276-281.

\bibitem{WGUO} Z.W. Wang, J.-M. Guo, A sharp upper bound of the nullity of a connected graph in terms of order and maximum degree, Linear Algebra Appl. 584(2020) 287-293.



\bibitem{wlt} Q. Wu, Y. Lu, B.-S. Tam, On connected signed graphs with rank equal to girth, Linear Algebra Appl. 651(2022) 90-115.

\bibitem{xzw} F. Xu, Q. Zhou, D. Wong, F.L. Tian, Complex unit gain graphs of rank 2, Linear Algebra Appl. 597(2020) 155-169.

\bibitem{yqt} G.H. Yu, H. Qu, J.H. Tu, Inertia of complex unit gain graphs, Appl. Math. Comput. 265(2015) 619-629.



\bibitem{ZF} F.Z. Zhang, Quaternions and matrices of quaternions, Linear Algebra Appl. 251(1997) 21--57.

\bibitem{QNZ} Q. N. Zhou, Y. Lu, Relation between the row left rank of a quaternion unit gain graph and the rank of its underlying graph, Electron. J. Linear Algebra. 39(2023) 181-198.



\bibitem{ZWS} Q. Zhou, D. Wong, D.Q. Sun, An upper bound of the nullity of a graph in terms of order and maximum degree, Linear Algebra Appl. 555(2018) 314-320.

\bibitem{zwt} Q. Zhou, D. Wong, B.-S. Tam, On connected graphs of order $n$ with girth $g$ and nullity $n-g$, Linear Algebra Appl. 630(2021) 56-68.



























\end{thebibliography}
\end{document}